\documentclass[smallextended]{svjour3_ax}       % onecolumn (second format)
\smartqed  % flush right qed marks, e.g. at end of proof
\usepackage{amsmath}
\usepackage{amssymb}
\usepackage{graphicx}

\newtheorem{thm}{Theorem}

\newtheorem{prop}{Proposition}

\newtheorem{rem}{Remark}

\newcommand{\Rset}{\mathbb{R}}

\newcommand{\calK}{\mathcal{K}}

\newcommand{\calP}{\mathcal{P}}
\newcommand{\calKL}{\mathcal{KL}}

\DeclareMathOperator*{\esssup}{ess\,sup}

\begin{document}

\title{Interpreting Models of Infectious Diseases in Terms of Integral Input-to-State Stability\thanks{
The work is supported in part by JSPS KAKENHI Grant Number 20K04536. }
}
%\subtitle{Do you have a subtitle?\\ If so, write it here}

\titlerunning{Integral Input-to-State Stability of Diseases Models}

\author{Hiroshi Ito}

\institute{Hiroshi Ito \at
              Department of Intelligent and Control Systems,  
              Kyushu Institute of Technology, \\ 
              680-4 Kawazu, Iizuka 820-8502, Japan \\
              \email{hiroshi@ces.kyutech.ac.jp}  %  \\
}

\date{Received: 2 May 2020 / Accepted: date}
\date{}

\maketitle

\begin{abstract}
This paper amis to develop a system-theoretic approach to 
ordinary differential equations which deterministically describe dynamics of 
prevalence of epidemics. The equations are treated as 
interconnections in which component systems are connected by signals.
The notions of integral input-to-state  stability (iISS) and 
input-to-state  stability (ISS) have been effective in addressing 
nonlinearities globally without domain restrictions in 
analysis and design of control systems. They provide 
useful tools of module-based methods integrating characteristics of component systems. 
This paper expresses fundamental properties of models of infectious diseases 
and vaccination through the language of iISS and ISS. 
The systematic treatment is expected to facilitate development of effective 
schemes of controlling the disease spread via non-conventional Lyapunov functions.

\end{abstract}

\section{Introduction}

For many decades, mathematical models of infectious diseases have been recognized as useful tools for 
public health decision-making during epidemics 
\cite{KermackEpidemics27,ANDMAYnature79,HETHinfdiseas,Keelinfdiseasbook09}. 
Detailed models help predict the future course of outbreak, while simple models allow one to 
understand 
mechanisms whose interpretations can lead to ideas of control strategies, such as 
vaccination, isolation, regulation and digital contact tracing or culling, 
which slow or ultimately eradicate the infection from the population. 
The objective of this paper is to facilitate the development of the latter. 
This paper does not report any novel behavior of disease transmission. 
Instead, this paper is devoted to a system and signal interpretation of 
behavior of classical and simple models 
of infectious diseases in the language of integral input-to-state stability (iISS) and 
input-to-state stability (ISS). It aims to take a first step toward development 
of an iISS/ISS-theoretic foundation for control design to combat infectious diseases. 
This paper reports that popular models share essentially the same qualitative 
behavior which can be analyzed and explained systematically via the same tools of iISS/ISS.

The notion of iISS and ISS have been accepted widely as mathematical tools 
to deal with and utilize nonlinearities effectively 
in the area of control \cite{SontagISS08,SONCOP,SontagSCL98}. 
The notions offer a systematic framework of module-based design of control systems. 
Once a system or a network is divided into `stable'' components, 
aggregating characteristics of components gives answers to control design problems 
systematically. 
The answers are global, and they are not restricted to small domains of variables. 
Without relying linearity, ISS allows one to handle systems 
based on boundedness of states with respect to bounded inputs. 
Importantly, the boundedness does not require finite operator gain, so that 
replacing linearity with ISS, we can handle a large 
class of nonlinear systems. However, nonlinearity such as bilinearity and saturation 
often prevent systems from being ISS. They are cases where nonlinearities retain 
convergence of state variables in the absence of inputs, but 
prevent state variables from being bounded in the presence of inputs. Such nonlinear systems 
are covered by iISS. 
Systems whose state is bounded for small inputs are grouped 
into the class of Strong iISS systems \cite{CHAANGITOStrISS14}. 
ISS and (Strong) iISS characterize both internal and external stability 
properties. The weak ``stability'' of (Strong) iISS 
components can be compensated by ISS components. This fact is one of 
useful and powerful tools of the iISS/ISS framework. 
Some of main ideas of iISS/ISS module-based arguments are 
packed in the iISS small-gain theorem \cite{ITOTAC06,ITOJIATAC09,ITOacc11neti,ANGSGiISS}
which is an extension of the 
ISS small-gain theorem \cite{JIA96LYA,KARJIAvectorSG09,LIUHILJIAcdc09,DRWsicon10,DASITOWIRejc11}. 
One of the important features of the small-gain methodology is 
that for interconnected systems and networks, it gives formula to explicitly construct 
non-conventional Lyapunov functions which not only establish stability properties of 
equilibria, but also properties with respect to external variables and parameters. 
For understanding behavior of diseases transmission, construction of Lyapunov functions has been 
one of major directions in mathematical epidemiology, 
and Lyapunov functions are known to be useful 
for analyzing global properties of stability of each given equilibrium 
(see \cite{KOROLyap02,KOROLyap04,FALLIDlypu07,OREGAN2010446,EnaNakIDlyapdelay11,NakEnaIDlyap14,SHUAIIDlyapu13,ChenLyapDI14} 
and references therein). 
However, the Lyapunov functions are classical, so that they 
do not unite characterizations of stability properties which vary with parameters and 
external elements. In other words, bifurcation analysis need to be 
performed separately to divide the Lyapunov function analysis into cases. 
In addition, it has been also common to make simplifying assumptions of 
precise population conservation that refuses 
external signals and perturbation of parameters 
(see, e.g., \cite{KOROLyap02,KOROLyap04,KOROgennonID06,OREGAN2010446,LIMULDseirlyapu95,EnaNakIDlyapdelay11}). 

This paper gives interpretations to behavior of typical models of infectious diseases 
in terms of basic characterizations of iISS and 
ISS \cite{SONISSV,SontagSCL98,ANGSONiISS} 
with the help of two tools provided by iISS and ISS. 
One tool is a criterion of the small-gain type. The other is a fusion of 
global and local nonlinear gain of components. 
The two tools formulated in this paper are not completely novel ideas, but
they extend standard concepts in the iISS/ISS methodology by 
specializing in the setup of diseases models. 
In addition to their usefulness, 
this paper illustrates how the same set of the basic characterizations and tools 
of iISS/ISS can be applied to popular models of infectious diseases uniformly
to help explain and capture their fundamental behavior globally without 
dividing the analysis into cases a prior. 

\section{Preliminaries}\label{sec:symbols}

Let the set of real numbers be denoted by $\Rset:=(-\infty,\infty)$. 
This paper uses the symbols $\Rset_+:=[0,\infty)$ and $\Rset_+^n:=[0,\infty)^n$. 
The convention $\infty\le\infty$ is used for notational simplicity. 
A function $\eta: \Rset_+\to\Rset_+$ is said to be of class $\calP$ and 
written as $\eta\in\calP$ if $\eta$ is continuous 
and satisfies $\eta(0)=0$ and $\eta(s)>0$ for 
all $s\in\Rset_+\setminus\{0\}$. 
A class ${\calP}$ function is said to be of 
class $\calK$ if it is strictly increasing, 
A class ${\calK}$ function is said to be 
of class ${\cal K}_{\infty}$ if it is unbounded. 
A continuous function $\beta: \Rset_+\times\Rset_+\to\Rset_+$ is 
said to be of class $\calKL$ if, for each fixed $t\geq 0$, 
$\beta(\cdot,t)$ is of class $\calK$ and, 
for each fixed $s>0$, $\beta(s,\cdot)$ is 
decreasing and $\lim_{t\to\infty}\beta(s,t)=0$. 
The zero function of appropriate dimension is denoted by $0$. 
Composition of $\eta_1, \eta_2: \Rset_+\to\Rset_+$ is expressed as 
$\eta_1\circ\eta_2$. 
For a continuous function $\eta: \Rset_+\to\Rset_+$, 
the function $\eta^\ominus$: 
$\Rset_+\to\overline{\Rset}_+:=[0,\infty]$ is defined as 
$\eta^\ominus(s)=\sup \{ \tau \in \Rset_+ : s \geq \eta(\tau) \}$. 
By definition, for $\eta\in\calK$, 
$\eta^\ominus(s)=\infty$ holds for all 
$s\geq \lim_{\tau\to\infty}\eta(\tau)$, and 
$\eta^\ominus(s)=\eta^{-1}(s)$ elsewhere. 
The set $\{1,2,3,\ldots,n\}$ is denoted by $\overline{n}$. 
For a set $U$, its cardinality is denoted by $|U|$. 
For $U\subset\overline{n}$ and $x\in\Rset^n$, 
the vector $x_U\in\Rset^{|U|}$ is the collection of $|U|$ components $x_i$, $i\in U$. 
%
%The symbol $\prod_{i=j}^k$ denotes products of sequences, and it is $1$ for $k< j$.

%
%%%%%%%%%%
%

A system of the form 
\begin{align}
\dot{x}(t)=f(x(t),u(t)) 
\label{eq:notationsys}
\end{align}
is said to be integral input-to-state stable (iISS) with respect to the input $u$ 
\cite{SontagSCL98} if 
there exist $\beta\in\calKL$, $\chi$, $\mu\in{\calK}\cup\{0\}$ such that, for all
measurable locally essentially bounded functions
$u:\Rset_+\to\Rset^p$, all $x(0)\in\Rset^n$ and all $t\ge 0$, its solution $x(t)$ exists 
and satisfies
\begin{align}
|x(t)|\le
\beta(|x(0)|, t) + \chi\left(\int_0^t\mu(|u(\tau)|)d\tau\right) ,  
\label{eq:defiISS}
\end{align}
where the symbol $|\cdot|$ denotes the Euclidean norm. 
System \eqref{eq:notationsys} is said to be 
strongly integral input-to-state stable (Strongly iISS) 
with respect to the input $u$ 
\cite{CHAANGITOStrISS14} if there exist $R>0$ and $\gamma\in{\calK}\cup\{0\}$ such that 
\begin{align}
&
{\esssup}_{t\in\Rset_+}|u(t)|<R\ \Rightarrow\ 
\nonumber\\
&\hspace{6ex}
|x(t)|\le
\beta(|x(0)|, t) + \gamma({\esssup}_{t\in\Rset_+}|u(t)|) 
\label{eq:defISS}
\end{align}
holds in addition to the above requirement \eqref{eq:defiISS}. 
The constant $R$ is called an input threshold. 
If the requirement \eqref{eq:defISS} is met with $R=\infty$, system \eqref{eq:notationsys} 
is said to be input-to-state stable (ISS) \cite{SONCOP}. 
The function $\gamma$ is called an ISS-gain function. 
An ISS system is Strongly iISS. A Strongly iISS system is iISS. 
Their converses do not hold true.  
iISS of \eqref{eq:notationsys} implies globally asymptotic stability of the equilibrium $x=0$ 
for $u=0$. 
If a radially unbounded and continuously differentiable function $V: \Rset^n\to\Rset_+$ 
satisfies 
\begin{align}
\forall x\in\Rset^n\hspace{1.5ex}
\forall u\in\Rset^p\hspace{1.5ex}
\frac{\partial V}{\partial x}(x)f(x,u)\le -\alpha(V(x))+\sigma(|u|) 
\label{eq:iISSLyap}
\end{align}
for some $\sigma\in\calK$ and some $\alpha\in\calP$ (resp., 
$\alpha\in\calK$ and $\lim_{s\to\infty}\alpha(s)\ge\lim_{s\to\infty}\sigma(s)$),
the function $V(x)$ is said to be an iISS (resp., ISS) Lyapunov function. 
The existence of an iISS (resp., ISS) Lyapunov function 
guarantees iISS (resp., ISS) of system \eqref{eq:notationsys} \cite{ANGSONiISS,SONISSV}. 
If an iISS Lyapunov function admits $\alpha\in\calK$ in \eqref{eq:iISSLyap}, the system is 
guaranteed to be Strong iISS \cite{CHAANGITOStrISS14}. 
iISS and ISS Lyapunov functions are conventional Lyapunov functions 
when the input $u$ is zero in system \eqref{eq:notationsys}. 
The Lyapunov-type characterization \eqref{eq:iISSLyap} yields 
\begin{align}
\limsup_{t\to\infty}V(x(t))\le \overline{\gamma}({\esssup}_{t\in\Rset_+}|u(t)|), 
\label{eq:asygain}
\end{align}
where $\overline{\gamma}:=\alpha^\ominus\circ\sigma$. 
If $u(t)\equiv u_0$ is a constant and \eqref{eq:iISSLyap} holds with an equality sign, 
we have $\lim_{t\to\infty}V(x(t))=\overline{\gamma}(|u_0|)$. 
The inequality \eqref{eq:asygain} is called the asymptotic gain property \cite{SONISSV}. 
If system \eqref{eq:notationsys} is not ISS, there exists 
no $\overline{\gamma}\in\calK$ satisfying 
the asymptotic gain property \eqref{eq:asygain}.
This is because \eqref{eq:defISS} never holds with $R=\infty$.
As in \eqref{eq:defiISS}, an iISS system which is not ISS accumulates 
its input, and the solution $x(t)$ exists for all $t\in\Rset_+$, but unbounded 
if the input is persistent. 
All the above are standard definitions given for sign-indefinite system \eqref{eq:notationsys}. 
When the vector field $f$ generates only non-negative $x(t)$ in 
\eqref{eq:notationsys} defined with $x(0)\in\Rset_+^n$ and $u(t)\in\Rset_+^p$, all the above 
definitions and facts are valid 
by replacing $\Rset$ with $\Rset_+$. 

\section{A Small-Gain Theorem for Generalized Balancing Kinetics}\label{sec:bilin}

This section shows a small-gain type method 
for establishing stability of dynamical networks in the framework of iISS. 
To propose a generalized formulation which includes 
a previously-developed criterion as a special case, consider 
$x(t)=[x_1(t),x_2(t),\ldots,x_n(t)]^T: \Rset_+\to\Rset_+^n$ governed by 
\begin{align}
\dot{x}_i=&
-\eta_{i-1,i}(x)+\sigma_{i,i-1}(x)
-\eta_{i+1,i}(x)
\nonumber \\[-.1ex]
& 
\hspace{2ex}
+\sigma_{i,i+1}(x) 
-\theta_i(x_i) 
+\kappa_i(w_i), \hspace{4ex}
i\in\overline{n}
\label{eq:sysbilin}
\end{align}
for any $x(0)\in\Rset_+^n$ and any 
measurable and locally essentially bounded function 
$w(t)=[w_1(t),w_2(t),\ldots,w_n(t)]^T\in\Rset_+^n$. 
In \eqref{eq:sysbilin}, the subscripts of $\eta$ and $\sigma$ 
are integers which are circular of length $n$. 
If a subscript $k$ by itself does not belong to $\overline{n}$, 
it stands for $((k-1)\,\mathrm{mod}\,n)+1$, 
where $i\,\mathrm{mod}\,n$ denotes the non-negative reminder of 
the division of an integer $i$ by $n$.
All subscripts in this section are circular. 
Assume that the functions $\eta_{i-1,i}$, $\eta_{i+1,i}$, 
$\sigma_{i,i-1}$, $\sigma_{i,i+1}: \Rset_+^n\to\Rset_+$ and 
$\theta_i\in\calP$ are locally Lipschitz and satisfy\footnote{
Under the assumption \eqref{eq:ellconv}, 
the implication \eqref{eq:sysbilinposi} is necessary and sufficient 
for guaranteeing $x(t)\in\Rset_+^n$.}
\begin{align}
x_i=0 
\hspace{1.8ex}\Rightarrow\hspace{1.8ex}
\eta_{i-1,i}(0)=\eta_{i+1,i}(0)=0  
\label{eq:sysbilinposi}
\end{align}
for all $i\in\overline{n}$. 
For all $i\in\overline{n}$, $\kappa_i\in\calK\cup\{0\}$ is assumed. 
Network \eqref{eq:sysbilin} is made of balancing mechanisms between 
state variables. 
The component $x_i$ is consumed at the rate 
$-\eta_{i+1,i}(x)$, and the consumption leads to the production 
of the downstream component $x_{i+1}$ at the rate $+\sigma_{i+1,i}(x)$. 
In the same way, the consumption 
$-\eta_{i-1,i}(x)$ of $x_i$ produces 
the upstream component $x_{i-1}$ at the rate $+\sigma_{i-1,i}(x)$. 
The rates are allowed to be
functions of $x$ instead of the local variable $x_i$. 
The extra rate $\theta_i$ of consumption in either the upstream or the 
downstream direction is a function of the local variable $x_i$. 
It is also important that the balancing between neighbors forms 
not only cycles of length $1$, but also cycles of length $n$. 
Assume that the rate functions in \eqref{eq:sysbilin} satisfy 
\begin{align}
\forall i\in\overline{n}
\hspace{1ex}
\forall j\in\{i-1,i+1\}
\hspace{1.5ex}
\exists \ell_{i,j}\ge 0 
\hspace{1.5ex}
\forall x\in\Rset_+^n
\hspace{1.5ex}
\sigma_{i,j}(x)\le\ell_{i,j}\eta_{i,j}(x). 
\label{eq:ellconv}
\end{align}
The following theorem can be proved. 

\begin{thm}\label{thm:sumBilin}
Assume that 
\begin{align}
&
\forall i\in\overline{n}
\hspace{1.5ex}
\ell_{i,i+1}\ell_{i+1,i} \le 1
\label{eq:ellsgcyc}
\\
&
\exists k\in\overline{n}
\hspace{1.5ex}
\ell_{k,k+1}\ell_{k+1,k}\le 
\prod_{i=1}^n
\dfrac{\ell_{i,i+1}}{\ell_{i+1,i}}
\le \frac{1}{\ell_{k,k+1}\ell_{k+1,k}} 
\label{eq:ellsgcycn}
\end{align}
are satisfied. 
Then network \eqref{eq:sysbilin} is iISS with respect the input $w$.  
If $\theta_i\in\calK$ (resp., $\theta_i\in\calK_\infty$) holds for all $i\in\overline{n}$ 
in addition, network \eqref{eq:sysbilin} is Strongly iISS 
(resp., ISS) with respect the input $w$. 
Furthermore, an iISS/ISS Lyapunov function $V(x)$ is 
\begin{align}
&
V(x)=\sum_{i\in\overline{n}}\lambda_ix_i
\label{eq:vsum}
\\
&
\forall i\!\in\!\{k\!+\!1,k\!+\!2,\ldots,k\!+\!n\}
\hspace{1.5ex}
\lambda_i=\!\!\prod_{j=k+2}^{i}\!\!
\sqrt{\frac{\ell_{j-1,j}}{\ell_{j,j-1}}} . 
\label{eq:lambdaBKlin}
\end{align}
\end{thm}

{\it Proof: }
Due to \eqref{eq:ellsgcyc}, 
from \eqref{eq:ellconv} and \eqref{eq:lambdaBKlin} it follows that 
\begin{align*}
\lambda_{i+1}\sigma_{i+1,i}(x) 
&=
\sqrt{\dfrac{\ell_{i,i+1}}{\ell_{i+1,i}}}\lambda_i
\sigma_{i+1,i}(x)
\\
&\le
\sqrt{\ell_{i,i+1}\ell_{i+1,i}}\lambda_i
\eta_{i+1,i}(x_{i})\
\\
&\le
\lambda_i\eta_{i+1,i}(x)  
\\
\lambda_{i}\sigma_{i,i+1}(x) 
&=
\sqrt{\dfrac{\ell_{i+1,i}}{\ell_{i,i+1}}}\lambda_{i+1}
\sigma_{i,i+1}(x)
\\
&\le
\lambda_{i+1}\eta_{i,i+1}(x) 
\end{align*}
for all $x\in\Rset_+^n$ at each 
$i\in\{k+1,k+2,\ldots,k+n-1\}$. 
By virtue of the second inequality in \eqref{eq:ellsgcycn}, 
properties \eqref{eq:ellconv} and \eqref{eq:lambdaBKlin} yield 
\begin{align*}
\lambda_{k}\sigma_{k,k+1}(x) 
&=
\prod_{i=k+1}^{k+n-1}
\sqrt{\dfrac{\ell_{i,i+1}}{\ell_{i+1,i}}}
\lambda_{k+1}\sigma_{k,k+1}(x)
\nonumber\\
&\le
\sqrt{\dfrac{\ell_{k+1,k}}{\ell_{k,k+1}}
\cdot\frac{1}{\ell_{k,k+1}\ell_{k+1,k}}
}\lambda_{k+1}
\sigma_{k,k+1}(x)
\\
&\le
\lambda_{k+1}\eta_{k,k+1}(x) .
\end{align*}
In the same way, 
the first inequality in \eqref{eq:ellsgcycn} gives 
\begin{align*}
\lambda_{k+1}\sigma_{k+1,k}(x) 
\le
\lambda_{k}\eta_{k+1,k}(x) . 
\end{align*}
Therefore, the function $V$ given in \eqref{eq:vsum} satisfies 
\begin{align*}
\dot{V}\le
\sum_{i\in\overline{n}}
\lambda_i\left(-\theta_i(x_i)+k_i(w_i)\right)
\end{align*}
along the trajectory $x(t)$ of network \eqref{eq:sysbilin}. Defining 
\begin{align*}
&
\alpha(s):=
\min_{\{x\in\Rset_+^n: V(x)=s\}}
\sum_{i\in\overline{n}}\theta_i(x_i)
\\
&
\sigma(s):=\sum_{i\in\overline{n}}\kappa_i(s)
\end{align*}
satisfies $\alpha\in\calP$, $\kappa_i\in{\calK}\cup\{0\}$ and 
\eqref{eq:iISSLyap}. Therefore, 
the function $V$ is an iISS Lyapunov function 
establishing iISS of network \eqref{eq:sysbilin}. 
We have $\alpha\in\calK$ establishing Strong iISS of network \eqref{eq:sysbilin} 
if $\theta_i\in\calK$ for all $i\in\overline{n}$. 
In the case of $\theta_i\in\calK_\infty$ for all $i\in\overline{n}$, 
we have $\alpha\in\calK_\infty$, and the function $V$ is 
an ISS Lyapunov function. 
\hfill $\square$

Recall that circulating subscripts of length $n$ are used for $\ell_{i,j}$ and $\lambda_i$. 
The above theorem extends a development of \cite{HITOconservcdc20} in which 
$\eta_{j,i}$ and $\sigma_{j,i}$ were restricted to be functions of $x_i$. 
The restriction disallows bilinearities and multiplicative nonlinearities to appear in 
$\eta_{j,i}$ and $\sigma_{j,i}$. Removal of the restrictions by 
Theorem \ref{thm:sumBilin} is the theoretical key in this paper. 

\begin{rem}\label{rem:noouterloop}
If there exists $k\in\overline{n}$ such that $\ell_{k,k+1}=0$ achieving 
\eqref{eq:ellconv}, the first inequality in \eqref{eq:ellsgcycn} is 
satisfied automatically with a sufficiently small $\ell_{k+1,k}>0$. 
The same applies to the case of $\ell_{k,k+1}=0$ with 
the second inequality in \eqref{eq:ellsgcycn}. 
\end{rem}
\begin{rem}\label{rem:sg}
Theorem \ref{thm:sumBilin} reduces to a special case of the iISS small-gain theorem 
for networks proposed in \cite{ITOacc11neti}. 
When the functions in \eqref{eq:sysbilin} are restricted to 
\begin{align}
&
\exists b_{i-i,i},b_{i+i,i}>0\hspace{1.5ex}
\forall x\in\Rset_+^n\hspace{1.5ex}
\nonumber\\
&\hspace{6ex}
b_{i-i,i}\eta_{i-i,i}(x)=b_{i+i,i}\eta_{i+i,i}(x)=\theta_i(x_i)
\\
&
\theta_i\in\calK
\\
&
\forall x\in\Rset_+^n\hspace{1.5ex}
\forall \hat{x}\in\Rset_+^n\cap\{\hat{x}_i=x_i\}\hspace{1.5ex}
\nonumber\\
&\hspace{6ex}
\sigma_{i,i-i}(x)\!=\!\sigma_{i,i-i}(\hat{x})
,\ 
\sigma_{i,i+i}(x)\!=\!\sigma_{i,i+i}(\hat{x})
\end{align}
for all $i\in\overline{n}$, network \eqref{eq:sysbilin} fit in the setup
of \cite{ITOacc11neti}, and 
the conditions 
\eqref{eq:ellsgcyc} and \eqref{eq:ellsgcycn} coincide with the cyclic small-gain 
condition presented in \cite{ITOacc11neti}, 
If all the functions in \eqref{eq:sysbilin} are non-zero, the network 
consists of $n$ simple directed cycles of length $1$, and $2$
simple directed cycles of length $n$. 
the inequality \eqref{eq:ellsgcyc} is the small-gain requirement for 
the cycles of length $1$, while 
the two inequalities in \eqref{eq:ellsgcycn} are for length $n$. 
\end{rem}

\section{Convergence via Zero Local ISS-Gain}\label{sec:gentheo}

This section focuses on an extra property which iISS systems often possess due to 
bilinear or multiplicative nonlinearities. To formulate it, 
consider $x(t)\in\Rset_+$ governed by 
\begin{align}
\dot{x}=f(x,w):=
\left[\begin{matrix}
f_1(x,w)\\
f_2(x,w)\\
\vdots\\
f_n(x,w)
\end{matrix}
\label{eq:sys}\right]
\end{align}
for any $x(0)\in\Rset_+^n$ and any 
measurable and locally essentially bounded function
$w(t)=[w_1(t),w_2(t),\ldots,w_p(t)]^T\in\Rset_+^p$. 
It is assumed that $f: \Rset_+^n\times\Rset_+^p\to\Rset$ is 
locally Lipschitz functions satisfying $f(0,0)=0$ and 
\begin{align}
x_i=0 
\ \Rightarrow\ 
\forall x\in\Rset_+^n\hspace{1.5ex}
\forall w\in\Rset_+^p\hspace{1.5ex}
f_i(x,w)\ge 0. 
\label{eq:sysposi}
\end{align}
Property \eqref{eq:sysposi} is necessary and sufficient for guaranteeing
$x(t)\in\Rset_+^n$ with respect to all $x(0)$ and $w$. 
The following proposition describes a property of iISS systems which 
completely reject the effect of small inputs on some state variables. 

\begin{prop}\label{prop:zerogain}
Suppose that system \eqref{eq:sys} is iISS with respect the input $w$, and 
admits a non-empty set $U\subset\overline{n}$, 
a radially unbounded and continuously differentiable function $V_U: \Rset^{|U|}\to\Rset_+$, 
and class $\calK$ functions $\alpha_U$, $\sigma_U$ such hat 
\begin{align}
\forall x\in\Rset_+^n\hspace{1.5ex}
\forall w\in\Rset_+^p\hspace{1.5ex}
\hspace{1.5ex}
\frac{\partial V_U}{\partial x_U}(x_U)f_U(x,w)\le -\alpha_U(V_U(x_U))+\sigma_U(w) . 
\label{eq:propxUISS}
\end{align}
If there exist a non-empty set $L\subset\overline{n}$, a real number $H>0$, and 
a radially unbounded and continuously differentiable function $V_L: \Rset^{|L|}\to\Rset_+$ 
such that for each $k\in(0,1)$, 
\begin{align}
\forall x\in\Rset_+^n\cap\left\{V_U(x_U)\le kH\right\}\hspace{1.5ex}
\forall w\in\Rset_+^p
\hspace{1.5ex}
\frac{\partial V_L}{\partial x_L}(x_L)f_L(x,w)\le -\alpha_L(V_L(x_L))
\label{eq:propxLdissip}
\end{align}
holds for some $\alpha_L\in\calP$, 
then the solution $x(t)$ of \eqref{eq:sys} satisfies 
\begin{align}
\sup_{t\in\Rset_+}w(t)< Q
\ \Rightarrow\  
\forall x(0)\in\Rset_+^n\hspace{1.5ex} 
\lim_{t\to\infty}x_L(t)=0  
\label{eq:propxLconv}
\end{align}
for $Q=\lim_{s\to H-}\sigma_U^\ominus\circ\alpha_U(s)$. 
\end{prop}

{\it Proof: }
Since system \eqref{eq:sys} is iISS, a solution 
$x(t)<\infty$ exists and unique for all $t\in\Rset_+$. 
Assume that a real number $H\in(0,\infty)$ satisfies 
\eqref{eq:propxLdissip}. 
Suppose that $\sup_{t\in\Rset_+}w(t)< Q$ for 
$Q=\lim_{s\to H-}\sigma_U^\ominus\circ\alpha_U(s)$. 
In the case of $Q=\infty$,  for each $w$, 
assumption \eqref{eq:propxUISS} implies the 
existence of $T\in\Rset_+$ and $\epsilon\in(0,H)$ 
such that
\begin{align}
\forall t\in[T,\infty)\hspace{1.5ex}
V_U(x_U(t))\le H-\epsilon . 
\label{eq:propultimatexU}
\end{align}
In the case of $Q<\infty$, since for each $w$, there exists $\delta>0$ satisfying 
$\sup_{t\in\Rset_+}w(t)\le Q-\delta$ , evaluating \eqref{eq:propxUISS} 
allows one to verify the existence of $T\in\Rset_+$ and $\epsilon\in(0,H)$ 
fulfilling \eqref{eq:propultimatexU} again. 
Therefore, the existence of $\alpha_L\in\calP$ satisfying 
\eqref{eq:propxLdissip} for each $k\in(0,1)$ 
ensures \eqref{eq:propxLconv}. 
\hfill $\square$

The above proposition makes use of property \eqref{eq:propxLdissip} which 
is the existence of a partial system completely rejecting the effect of 
its interconnecting inputs whose magnitude is 
smaller than a threshold. 
Let $L^C=\overline{n}\setminus L$. 
For $x_L$-system defined by $\dot{x}_L=f_L(x,w)$,
$x_{L^C}(t)$ and $w(t)$ are exogenous. 
Property \eqref{eq:propxLdissip} implies 
that the ISS-gain from the input $[x_{L^C}^T, w^T]^T$ to the state $x_L$ is zero as long 
as $V_U(x_U)<H$. Hence, in Proposition \ref{prop:zerogain}, 
$x_L$-system is required to admit zero 
ISS-gain locally, although it is not required to be ISS globally. 
In the case of $L=\overline{n}$, property \eqref{eq:propxLdissip} implies 
Strong iISS of $x_L$-system since iISS is assumed in Proposition \ref{prop:zerogain}.
This property \eqref{eq:propxLdissip} is very common, although it may not have 
been focused very often in the literature of systematic control methodology. 
For example, any of scalar 
non-negative systems
\begin{align*}
&
\dot{\xi}=-\xi+\xi w
\\
&
\dot{\xi}=-\frac{\xi}{1+\xi}+\frac{\xi w}{1+\xi}  
\\
&
\dot{\xi}=-\frac{\xi}{1+\xi^2}+\frac{\xi w}{1+\xi^2}  
\end{align*}
has the threshold at $H=1$ for $x_L=\xi$ and $V(x_L)=x_L$. 
The state $\xi$ converges to zero for $w<1$, while 
$\xi$ increases as long as $w>1$, although the increase is 
within the property of iISS . 
The above three systems are Strong iISS. The zero local ISS-gain they actually 
have is ``stronger'' than Strong iISS. 
In this way, the threshold is a bifurcation point bringing in a superior
stability property. 
Bilinearity and multiplicative nonlinearity give rise to the bifurcation, and 
Proposition \ref{prop:zerogain} aims at highlighting such systems.

\begin{rem}\label{rem:estimates}
Since Proposition \ref{prop:zerogain} assumes the entire system \eqref{eq:sys} to be 
only iISS, the variable $x_L(t)$ can increase until $V_U(x_U(t))<H$ 
is achieved. Then an outbreak of $x_L$ occurs, and $x_L(t)$ exhibits a peak. 
Property \eqref{eq:propxUISS} allows one to estimate that 
the smaller $w$ is, the shorter the time when $x_L$ starts decreasing 
(the duration of the growth phase) is. 
The larger $Q$ or $H$ is, the shorter the growth phase is. 
The increasing rate of $x_L(t)$ until $V_U(x_U(t))<H$ can be 
estimated by the iISS property of the entire system, or   
$x_L$-system, which is more specific than the entire system. 
For example, an upper bound of 
$({\partial V_L}/{\partial x_L})f_L$ for $V_U(x_U)> H$ gives such 
information, which is not stated explicitly in 
Proposition \ref{prop:zerogain}. 
\end{rem}

\section{SIR Model}\label{sec:sirv}

Consider the solution $x(t):=[S(t), I(t), R(t)]^T\in\Rset_+^3$ of the 
ordinary differential equation 
\begin{subequations}\label{eq:sirv}
\begin{align}
\dot{S}=&B -\mu S -\beta IS
\label{eq:sirvS}
\\
\dot{I}=&\beta IS -\gamma I-\mu I
\label{eq:sirvI}
\\
\dot{R}=&\gamma I-\mu R
\label{eq:sirvR}
\end{align}
\end{subequations}
defined for any $[S(0), I(0), R(0)]^T\in\Rset_+^3$ and 
any measurable and locally essentially 
bounded function $B: \Rset_+\to\Rset_+$ . 
The variable $S(t)$ describes the (continuum) number of susceptible population. 
$I(t)$ is the number of infected individuals, while 
$R(t)$ is of individuals recovered with immunity. 
$B(t)$ is the newborn/immigration rate.  
The positive number $\beta$, $\gamma$ and $\mu$ are parameters describing 
the contact rate, the recovery rate and the death rate, respectively. 
The equation \eqref{eq:sirv} is 
referred to as the (classic) SIR endemic model \cite{{HETHinfdiseas}}. 
When $B=0$ and $\mu=0$, the equation is referred to as the (classic) SIR epidemic model. 

S-system \eqref{eq:sirvS} and R-system \eqref{eq:sirvR} are ISS 
with respect to the input $[B, I]^T$ and the input $I$, respectively. 
In fact, the positive variables $S$ and $R$ by themselves are ISS Lyapunov functions, 
due to the definition \eqref{eq:iISSLyap}. 
The ISS property is also clear since R-system \eqref{eq:sirvR} is linear, and 
S-system \eqref{eq:sirvS} is bounded from above by the solution of 
the linear system $\dot{S}=B -\mu S$. The two linear systems are 
asymptotically stable. 
By contrast, due to the bilinear term $\beta IS$, I-system \eqref{eq:sirvI} is not 
ISS with respect to the input $S$ since \eqref{eq:sirvI} generates 
unbounded $I(t)$ for any constant input $S > (\gamma+\mu)/\beta$. 
I-system \eqref{eq:sirvI} is Strong iISS since $V_I(I)=\log(1+I)$ yields 
\begin{align*}
\dot{V}_I \le  -\frac{(\gamma+\mu) I}{1+I}+\beta S . 
\end{align*}
Hence, the SIR model \eqref{eq:sirv} is a cascade consisting of 
an ISS system, an Strong iISS system and an ISS system. 
The argument of iISS cascade analysis 
\cite{CHAIANGmtns06,ITOTAC06,ITOTAC10,CHAANGITOStrISScas14,HITOcasSICE17} 
can show that the SIR model \eqref{eq:sirv} is Strong iISS with respect to 
the input $B$. This stability assessment is not yet very informative, 
although it holds true. The rest this section demonstrates that 
the developments in Sections \ref{sec:bilin} and \ref{sec:gentheo} 
improve the stability assessment for capturing discriminative behavior of 
the spread of infectious diseases.

Theorem \ref{thm:sumBilin} establishes that the SIR model \eqref{eq:sirv} is ISS, 
which is ``stronger'' than Strong iISS, although the SIR model involves 
I-system which is not ISS. 
In fact,  the model \eqref{eq:sirv} satisfies \eqref{eq:sysbilin} and 
\eqref{eq:ellconv} with 
\begin{align}
\ell_{i+1,i}=\ell_{i,i+1}=1, \ \theta_i(s)=\mu s, 
\ i=1,2,3. 
\label{eq:ellmusierv}
\end{align}
These parameters satisfy \eqref{eq:ellsgcyc} and \eqref{eq:ellsgcycn}, and yield 
$\lambda_1=\lambda_2=\lambda_3=1$. 
Theorem \ref{thm:sumBilin} gives 
the total number $N(t):=S(t)+I(t)+R(t)=V(x(t))$ as an ISS Lyapunov function for 
the SIR model \eqref{eq:sirv}. Indeed, 
\begin{align}
\dot{N} =  -\mu N + B
\label{eq:NdissipISS}
\end{align}
is satisfied along the solution $x(t)$ of \eqref{eq:sirv}. 
The above inequality implies 
\begin{align}
\lim_{t\to\infty}N(t)= \frac{B}{\mu} 
\label{eq:NdissipISSUB}
\end{align}
for any constant $B\in\Rset_+$. 
If $B$ is not constant, we have the asymptotic gain property 
\begin{align}
\limsup_{t\to\infty}N(t)\le \frac{1}{\mu}\esssup_{t\in\Rset_+}B(t). 
\label{eq:NdissipISSUBt}
\end{align}

SI-system consisting of \eqref{eq:sirvS} and \eqref{eq:sirvI} 
is ISS with respect to the input $B$ since defining the sum 
$V_{U}(x_U(t))=S(t)+I(t)$ yields 
\begin{align*}
\dot{V}_{U} =  -\mu V_{U}-\gamma I + B  
\end{align*}
for $x_U=[S,I]^T$. 
This confirms \eqref{eq:propxUISS} with 
$\alpha_U(s)=\mu s$ and $\sigma_U(s)=s$. 
I-system defined by \eqref{eq:sirvI} fulfills \eqref{eq:propxLdissip} with 
$x_L=I$ and 
\begin{align}
H:=\frac{\gamma+\mu}{\beta}
\label{eq:defRsirv}
\end{align}
since $S \le V_U(x_U)$. 
Recall that  ISS of the entire model \eqref{eq:sirv} implies iISS of \eqref{eq:sirv}. 
Hence, Proposition \ref{prop:zerogain} with $Q=\mu H$ guarantees 
\begin{align}
\esssup_{t\in\Rset_+}{B}(t)
<\mu H
\ \Rightarrow\ \lim_{t\to\infty}I(t)=0 . 
\end{align}
The convergence of $I$ to zero implies 
$\lim_{t\to\infty}R(t)=0$ since \eqref{eq:sirvR} is ISS. Indeed, 
it is a stable linear system. Thus 
\begin{align}
\esssup_{t\in\Rset_+}{B}(t)
<\mu H
\ \Rightarrow\ \lim_{t\to\infty}I(t)=\lim_{t\to\infty}R(t)=0 . 
\end{align}
For any constant $B\in\Rset_+$, 
by virtue of \eqref{eq:NdissipISSUB} or \eqref{eq:sirvS}, we arrive at 
\begin{align}
{B}<\mu H
\ \Rightarrow\ \lim_{t\to\infty}x(t)=\left[\frac{B}{\mu},0,0\right]^T . 
\label{eq:sirv_epidemicequi}
\end{align}
As Proposition \ref{prop:zerogain} is applied in the above analysis, 
regarding $S$ as a parameter, I-system has a bifurcation point at $S=H$. 
The origin $I=0$ is an unstable equilibrium if $S>H$. 
As explained in Remark \ref{rem:estimates}, $I(t)$ increases until $S(t) \le H$. 
It is not a possibility any more and the growth phase occurs 
since I-system is not ISS. 
In the case of ${B}<\mu H$ and $S(0) > H$, if $I(0)>0$, the infected population  
$I(t)$ peaks before converging to zero. 
The smaller $B$ and the larger $(\gamma+\mu)/\beta$ are, the shorter the time 
to the peak is. 
The larger $\mu$ and the smaller $\beta$ are, the smaller the  
increase rate of $I$ in the growth phase is. 

\begin{figure}[t]
%\hspace{-1ex}%
\begin{center}
\includegraphics[width=9.2cm,height=5.4cm]{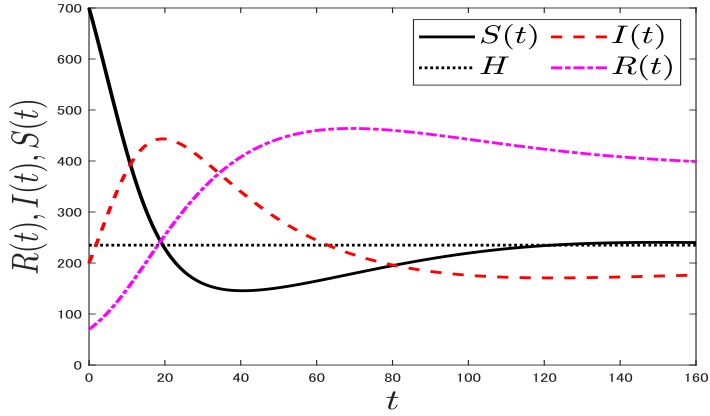}
\vspace{-1.2ex}
\end{center}
\caption{Populations of the SIR model \eqref{eq:sirv} with $B(t)\equiv 12$, 
which is $R_0=B/(\mu H)=3.4043>1$.}
\label{fig:sire}
%\vspace{2ex}
\end{figure}

\begin{figure}[t]
%\hspace{-1ex}%
\begin{center}
\includegraphics[width=9.2cm,height=5.4cm]{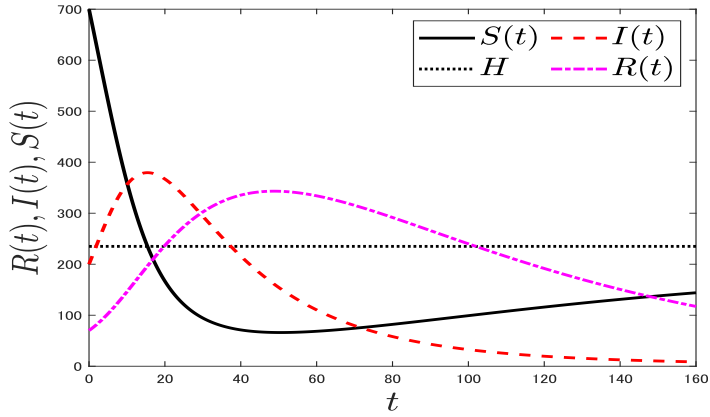}
\vspace{-1.2ex}
\end{center}
\caption{Populations of the SIR model \eqref{eq:sirv} with $B(t)\equiv 3$, 
which is $R_0=B/(\mu H)=0.8511<1$.}
\label{fig:sirf}
%\vspace{2ex}
\end{figure}

Typical responses of the SIR model \eqref{eq:sirv} are shown in 
Fig. \ref{fig:sire} and Fig. \ref{fig:sirf} for 
$\beta=0.0002$, $\mu=0.015$ and $\gamma=0.032$ with 
$S(0)=700$, $I(0)=200$ and $R(0)=70$. 
For a general newborn/immigration rate $B$, 
the SIR model \eqref{eq:sirv} is ISS. It means that 
the disease can remain as endemic. That is, the infected population 
$I(t)$ is bounded, but it can become very large, 
and $\liminf_{t\to\infty}I(t)>0$ can hold. 
Since $I$-system is Strong iISS and it is not ISS, 
the infected $I(t)$ starts with an growth phase unless 
the initial susceptible population is below the threshold $H=235$. 
The infected population $I(t)$ never decreases to zero if 
the newborn/immigration rate $B$ is above the threshold $\mu H=3.525$ (Fig. \ref{fig:sire}). 
If the newborn/immigration rate $B$ does not exceed the threshold, 
the convergence of $I(t)$ to zero is guaranteed, and 
the disease is eradicated (Fig. \ref{fig:sirf}). 
This bifurcation is the central feature of the SIR model \eqref{eq:sirv}. 

In mathematical epidemiology, for constant $B\in\Rset_+$, the value 
\begin{align}
R_0:=
\frac{\beta B}{\mu(\gamma+\mu)}=
\frac{B}{\mu H}
\end{align}
is called the basic reproduction number. 
Solving the simultaneous equation of $\dot{S}=0$, 
$\dot{I}=0$ and $\dot{R}=0$ in \eqref{eq:sirv} 
with constant $B\in\Rset_+$ gives the steady-state 
value in \eqref{eq:sirv_epidemicequi} and 
\begin{align}
&
R_0\ge 1
\ \Rightarrow\ 
x_e=\left[
H, 
\frac{\mu(R_0-1)}{\beta}, 
\frac{\gamma(R_0-1)}{\beta}
\right]^T.
\label{eq:sirv_endemicequi}
\end{align}
The steady-state value in \eqref{eq:sirv_epidemicequi} is called 
the disease-free equilibrium, while 
$x_e$ in \eqref{eq:sirv_endemicequi} is called 
the endemic equilibrium. 
Since $x(t)\in\Rset_+$, 
the endemic equilibrium exists only if $R_0\ge 1$. 
For $R_0=1$, the endemic equilibrium is identical with 
the disease-free equilibrium. 
The endemic equilibrium is consistent with \eqref{eq:NdissipISSUB} since 
the definition of $R_0$ yields 
\begin{align}
H+
\frac{\mu(R_0-1)}{\beta}+
\frac{\gamma(R_0-1)}{\beta}
=
\frac{(\gamma+\mu)R_0}{\beta}
=\frac{B}{\mu} . 
\end{align}
Recall that the entire SIR model \eqref{eq:sirv} and SI-systems are ISS. 
Since iISS of I-system 
which is not ISS accumulates the amount $S(t)-H$, 
there exists a unbounded increasing sequence $\{t_i\}_{i\in\{0,1,2,...\}}$ in 
$\Rset_+$ such that 
\begin{align}
\forall x(0)\in\Rset_+\hspace{-.3ex}\setminus\!\{I(0)>0\}\hspace{1.5ex}
\lim_{i\to\infty}x_2(t_i)=x_{e,2} 
\end{align}
if $R_0> 1$ for constant $B\in\Rset_+$,  where $x_e=[x_{e,1},x_{e,2},x_{e,3}]^T$. 

\begin{rem}
Many analytic studies on disease models have made assumptions 
on $B$ to make the analysis simple. 
For example, assuming $B=\mu(S+I+R)$ or 
an equivalent setup in the SIR model results in 
$S(t)+I(t)+R(t)=N$ for all $t$ with a positive constant $N$ 
\cite{KOROLyap02,KOROgennonID06,OREGAN2010446}. 
The mutual dependence between the variables removes one of the three variable 
from \eqref{eq:sirv}. Many studies on variants of the SIR model have also 
relied on this significant simplification 
(e.g., \cite{LIMULDseirlyapu95,KOROLyap04,EnaNakIDlyapdelay11}). 
This paper does not make such assumptions since we are 
interested in not only avoiding the limited usefulness of analysis, 
but also the perturbation of the newborn rate $B$. This is the reason 
why this paper explicitly refers to $B$ as the immigration rate, which is 
an exogenous signal, while the assumption $B=\mu(S+I+R)$ forces 
the variable $B$ to be endogenous completely. 
In fact, with the simplification, the robustness notion this paper pursues does not 
arise at all. 
\end{rem}

\section{SEIS Model}\label{sec:seis}

Consider $x(t):=[S(t), E(t), I(t)]^T\in\Rset_+^3$ governed by 
\begin{subequations}\label{eq:seis}
\begin{align}
\dot{S}=&B -\mu S -\beta IS +\gamma I
\label{eq:seisS}
\\
\dot{E}=&\beta IS -\epsilon E-\mu E
\label{eq:seisE}
\\
\dot{I}=& \epsilon E-\gamma I-\mu I
\label{eq:seisI}
\end{align}
\end{subequations}
with $[S(0), E(0), I(0)]^T\in\Rset_+^3$ and 
$B: \Rset_+\to\Rset_+$. 
The equation \eqref{eq:seis} is 
referred to as the SEIS model \cite{KOROLyap04}. 
The SEIS model is known to be useful for describing diseases which 
have non-negligible incubation periods. 
The variable $E$ represents the (continuum) number of infected individuals 
who are not yet infectious. 
The SEIS model also consider infections which do not give 
long lasting immunity, and recovered individuals become
susceptible again.
As seen in \eqref{eq:seis}, 
the short immunity forms a circle of length\footnote{A 
notion in graph theory} $3$, which the SIR model 
does not have. 

The SEIS model \eqref{eq:seis} satisfies \eqref{eq:sysbilin} and 
\eqref{eq:ellconv} with \eqref{eq:ellmusierv}. Conditions 
\eqref{eq:ellsgcyc} and \eqref{eq:ellsgcycn} are satisfied. Thus 
Theorem \ref{thm:sumBilin} establishes ISS of the SEIS model \eqref{eq:seis}
with respect to the input $B$, and an ISS Lyapunov function is obtained as 
\eqref{eq:vsum} with $\lambda_1=\lambda_2=\lambda_3=1$, i.e., 
$V(x)=S+E+I=:N$. In fact, properties \eqref{eq:NdissipISS}, 
\eqref{eq:NdissipISSUB} and \eqref{eq:NdissipISSUBt} are verified. 

\begin{figure}[t]
%\hspace{-1ex}%
\begin{center}
\includegraphics[width=9.2cm,height=5.4cm]{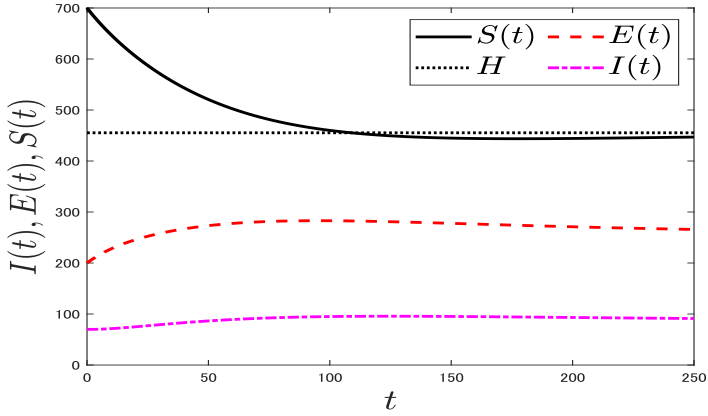}
\vspace{-1.2ex}
\end{center}
\caption{Populations of the SEIS model \eqref{eq:seis} with $B(t)\equiv 12$, 
which is $R_0=B/(\mu H)=1.757>1$.}
\label{fig:seise}
%\vspace{2ex}
\end{figure}

\begin{figure}[t]
%\hspace{-1ex}%
\begin{center}
\includegraphics[width=9.2cm,height=5.4cm]{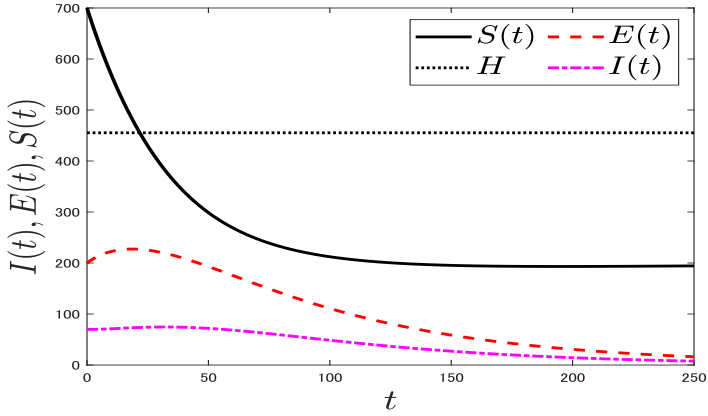}
\vspace{-1.2ex}
\end{center}
\caption{Populations of the SEIS model \eqref{eq:seis} with $B(t)\equiv 3$, 
which is $R_0=B/(\mu H)=0.4393<1$.}
\label{fig:seisf}
%\vspace{2ex}
\end{figure}

For $x_U:=x=[S, E, I]^T$ and $V_U=V$, 
property \eqref{eq:NdissipISS} achieves 
\eqref{eq:propxUISS} with 
$\alpha_U(s)=\mu s$ and $\sigma_U(s)=s$. 
EI-system consisting of \eqref{eq:seisE} and \eqref{eq:seisI} is iISS 
with respect to the input $S$ since 
the choice $V_{EI}=\log(1+E+I)$ satisfies 
\begin{align}
\dot{V}_{EI} \le  -\frac{\mu (E+I)}{1+E+I} + \beta S 
\label{eq:EIdissipiISS}
\end{align}
along the solution $[E(t), I(t)]^T$ of EI-system.
EI-system is, however, not ISS. EI-system is a Strongly iISS system 
admitting a zero local ISS-gain. 
To see this, one can make use of the Lyapunov function proposed in Theorem \ref{thm:sumBilin} 
by regarding $S$ as a constant $S^\sharp$ for EI-system. 
EI-system satisfies \eqref{eq:sysbilin} and \eqref{eq:ellconv} with 
\begin{align}
& 
\ell_{1,2}=\frac{\beta S^\sharp}{a(\gamma+\mu)}, \quad 
\ell_{2,1}=\frac{\epsilon}{a(\epsilon+\mu)}
\\
& 
\theta_1(s)=(1\!-\!a)(\epsilon\!+\!\mu)s
, \ 
\theta_2(s)=(1\!-\!a)(\gamma\!+\!\mu)s
\label{eq:ellseisei}
\end{align}
for an arbitrarily given $a\in(0,1)$. Defining $V_{L}(x_L)$ for 
$x_L=[E, I]^T$ as in \eqref{eq:vsum} and \eqref{eq:lambdaBKlin} gives 
\begin{align}
V_{L}(x_L)=E+\lambda_2 I
, \quad 
\lambda_2=
\sqrt{\frac{\beta S^\sharp(\epsilon+\mu)}
{(\gamma+\mu)\epsilon}} . 
\label{eq:defEWEI}
\end{align}
Then along the solution of \eqref{eq:seisE} and \eqref{eq:seisI} we obtain 
\begin{align}
\frac{d}{dt}V_L(x_L)=&
\sqrt{\frac{\beta\epsilon(\epsilon+\mu)}{\gamma+\mu}}
\left(
\sqrt{S^\sharp}-
\sqrt{\frac{(\gamma+\mu)(\epsilon+\mu)}{\beta\epsilon}}
\right)E
\nonumber\\
&
+
\sqrt{\beta^2 S^\sharp}
\left(
\frac{S}{\sqrt{S^\sharp}}-
\sqrt{\frac{(\gamma+\mu)(\epsilon+\mu)}{\beta\epsilon}}
\right)I . 
\label{eq:EidissipISS}
\end{align}
Define 
\begin{align}
H:=\frac{(\gamma+\mu)(\epsilon+\mu)}{\beta\epsilon} . 
\label{eq:seisboint}
\end{align}
Equation \eqref{eq:EidissipISS} allows one to see that $S=H$ 
is a bifurcation point of EI-system. To this end, 
pick $S^\sharp=a^2H$. 
Due to $S \le V_U(x_U)$, 
EI-system satisfies \eqref{eq:propxLdissip} with 
\eqref{eq:seisboint} 
since\footnote{This is consistent with \eqref{eq:ellsgcyc} and 
\eqref{eq:ellsgcycn} evaluated with the existence of $a\in(0,1)$ in 
Theorem \ref{thm:sumBilin}.} 
such a parameter $a\in(0,1)$ can be taken for each $k\in(0,1)$. 
On the other hand, if $S$ is a constant satisfying $S>H$, 
equation \eqref{eq:EidissipISS} with $S^\sharp=S$ yields 
${d}V_L(x_L)/{dt}>0$ unless $E=I=0$. 
Therefore, EI-system is not ISS, but the convergence property \eqref{eq:propxLdissip} 
is met. Hence,
Proposition \ref{prop:zerogain} can be invoked for $Q=\mu H$, 
and concludes that 
\begin{align}
\sup_{s\in\Rset_+}{B}(t)<\mu H
\ \Rightarrow\ \lim_{t\to\infty}\!E(t)=\lim_{t\to\infty}\!I(t)=0. 
\end{align}
From \eqref{eq:NdissipISSUB}, 
the solution $x(t)=[S(t), E(t), I(t)]^T$ of \eqref{eq:seis}
satisfies \eqref{eq:sirv_epidemicequi} for any constant $B$. 
Accordingly, all the observations made for the SIR model apply to the SEIS model except that 
the increase of the infectious individuals $I$ is replaced by the increase of 
the weighted sum of infected individuals $E+\lambda_2 I$. 
Simulations of the SEIS model \eqref{eq:seis} are shown in 
Figs. \ref{fig:seise} and \ref{fig:seisf} for 
$\beta=0.0002$, $\mu=0.015$ and $\gamma=0.032$ with 
$S(0)=700$, $I(0)=200$ and $R(0)=70$. The thresholds are $H=455.3$ 
and $\mu H=6.83$. The disease is removed in 
Fig. \ref{fig:seisf} plotted for $B=3<\mu H$. 
The infected and infectious populations remain high 
in Fig. \ref{fig:seise} computed for $B=12>\mu H$. 

\section{MSIR and SEIR Models}\label{sec:msir}

Consider $x(t):=[M(t), S(t), I(t), R(t)]^T\in\Rset_+^4$ for 
\begin{subequations}\label{eq:msir}
\begin{align}
\dot{M}=&B -\delta M -\mu M
\label{eq:msirM}
\\
\dot{S}=&\delta M -\mu S -\beta IS
\label{eq:msirS}
\\
\dot{I}=&\beta IS -\gamma I-\mu I
\label{eq:msirI}
\\
\dot{R}=&\gamma I-\mu R
\label{eq:msirR}
\end{align}
\end{subequations}
which is called the MSIR model \cite{HETHinfdiseas,MCLANDmsirID88}. 
The variable $M$ represents delay in 
becoming susceptible due to the maternally derived immunity. 
The analysis of the MSIR model is almost the same as that of the SIR model. 
With 
\begin{align}
\ell_{i+1,i}=\ell_{i,i+1}=1, \ \theta_i(s)=\mu s, 
\ i=1,2,3,4, 
\label{eq:ellmumsir}
\end{align}
Theorem \ref{thm:sumBilin} assures that 
the function $N(t):=M(t)+S(t)+I(t)+R(t)$ proves ISS of 
\eqref{eq:msir}, and the ultimate bounds \eqref{eq:NdissipISSUB} and 
\eqref{eq:NdissipISSUBt} via 
\eqref{eq:NdissipISS}. 
Thus, the choices $V_U(x)=N$ and $x=x_U$ give \eqref{eq:propxUISS} with 
$\alpha_U(s)=\mu s$ and $\sigma_U(s)=s$. 
Because of the bilinear term $\beta IS$, 
I-system in \eqref{eq:msir} is not ISS, but Strongly iISS. 
Thus, the variable $I(t)$ increases until $S(t) \le H$, where 
the bifurcation point $H$ is defined as \eqref{eq:defRsirv}. 
The ISS gain of I-system is zero for the input $\sup_{s\in\Rset_+}S(t)<H$. 
In fact, I-system satisfies \eqref{eq:propxLdissip}. 
Thus, Proposition \ref{prop:zerogain} establishes 
\begin{align}
{B}<\mu H
\ \Rightarrow\ \lim_{t\to\infty}x(t)=\left[\frac{B}{\mu},0,0,0\right]^T 
\label{eq:sirv_epidemicequi4}
\end{align}
for the MSIR model \eqref{eq:msir} in the case of constant $B\in\Rset_+$. 

%%%%%%%%%%

Finally, the SEIR model consists of 
\begin{subequations}\label{eq:seir}
\begin{align}
\dot{S}=&B -\mu S -\beta IS
\label{eq:seirS}
\\
\dot{E}=&\beta IS -\epsilon E-\mu E
\label{eq:seirE}
\\
\dot{I}=& \epsilon E-\gamma I-\mu I
\label{eq:seirI}
\\
\dot{R}=& \gamma I-\mu R . 
\label{eq:seirR}
\end{align}
\end{subequations}
Its state vector is $x(t):=[S(t), E(t), I(t), R(t)]^T\in\Rset_+^4$ \cite{LIMULDseirlyapu95}. 
The SEIR model can be analyzed as done for the SIES model. 
Theorem \ref{thm:sumBilin} qualifies 
$N(t):=S(t)+E(t)+I(t)+R(t)$ as an ISS Lyapunov function proving 
ISS of \eqref{eq:seir}, and provides 
the ultimate bounds \eqref{eq:NdissipISSUB} and \eqref{eq:NdissipISSUBt} via 
\eqref{eq:NdissipISS}. 
Property \eqref{eq:propxUISS} holds with 
$\alpha_U(s)=\mu s$ and $\sigma_U(s)=s$ for $V_U(x)=N$ and $x=x_U$. 
EIR-system consisting of \eqref{eq:seirE}, \eqref{eq:seirI} and \eqref{eq:seirR} 
is Strongly iISS with respect to the input $S$, although 
the bilinear term $\beta IS$ prevent it from being ISS. 
Using $V_L(x_L)=\lambda_E E+\lambda_I I+\lambda_R R$ 
for appropriate $\lambda_E$, $\lambda_I$, $\lambda_R>0$ 
given by Theorem \ref{thm:sumBilin}, 
one can  show that 
EIR-system is not ISS. 
The function $V_L(x_L)$ also shows that EIR-system 
admits the zero ISS gain for the input $\sup_{s\in\Rset_+}S(t)<H$, where 
the bifurcation point $H$ is defined as \eqref{eq:seisboint}. 
EIR-system satisfies \eqref{eq:propxLdissip}. 
Therefore, Proposition \ref{prop:zerogain} concludes that 
the SEIR model \eqref{eq:seir} satisfies \eqref{eq:sirv_epidemicequi4} 
for constant $B\in\Rset_+$. 

\section{Vaccination models}\label{sec:vacci}

One way of eradicating infectious diseases is to vaccinate newborns and entering individuals. 
Let the constant $P\in(0,1)$ denote the vaccination fraction. Considering 
a vaccine giving lifelong immunity \cite{HETHinfdiseas}, the SIR model can be modified as 
\begin{subequations}\label{eq:sirvvacci}
\begin{align}
\dot{S}=&B(1-P) -\mu S -\beta IS
\label{eq:sirvvacciS}
\\
\dot{I}=&\beta IS -\gamma I-\mu I
\label{eq:sirvvacciI}
\\
\dot{R}=&\gamma I-\mu R
\label{eq:sirvvacciR}
\\
\dot{A}=&BP-\mu A, 
\label{eq:sirvvacciA}
\end{align}
\end{subequations}
where $A$ is the number of vaccinated individuals. 
Since \eqref{eq:ellconv} is satisfied with \eqref{eq:ellmumsir}, 
Theorem \ref{thm:sumBilin} assures that 
$N(t):=S(t)+I(t)+R(t)+A(t)$ is an ISS Lyapunov function for 
the model \eqref{eq:sirvvacci}, 
and establishes the ultimate bounds \eqref{eq:NdissipISSUB} and 
\eqref{eq:NdissipISSUBt} via 
\eqref{eq:NdissipISS}. 
S-system is ISS with respect to the input $B$. In fact, 
property \eqref{eq:propxUISS} is met with $V_U(x_U)=S$, $x_U=S$, 
$\alpha_U(s)=\mu s$ and $\sigma_U(s)=(1-P)s$. 
IR system is the same as that of the SIR model. Hence, a 
bifurcation point $H$ is obtained as \eqref{eq:defRsirv}. 
The convergence of $I$ to zero implies 
$\lim_{t\to\infty}R(t)=0$. The variable $A$ reaches its steady state 
since $A$-system is ISS. Indeed, 
Hence, Proposition \ref{prop:zerogain} guarantees 
\begin{align}
{B}(1\!-\!P)<\mu H
\ \Rightarrow\ \lim_{t\to\infty}x(t)=\left[\frac{B(1\!-\!P)}{\mu},0,0,\frac{BP}{\mu}\right]^T 
\label{eq:sirv_epidemicequivacciA}
\end{align}
for any constant $B\in\Rset_+$ and $P\in(0,1)$. 
Hence, a vaccination fraction $P$ which is sufficiently close to $1$ 
can eradicate the disease.

%%%%%%%%%%

Another way to model the newborn vaccination within the SIR model is 
\begin{subequations}\label{eq:sirvvaccir}
\begin{align}
\dot{S}=&B(1-P) -\mu S -\beta IS
\label{eq:sirvvaccirS}
\\
\dot{I}=&\beta IS -\gamma I-\mu I
\label{eq:sirvvaccirI}
\\
\dot{R}=&\gamma I-\mu R +BP . 
\label{eq:sirvvaccirR}
\end{align}
\end{subequations}
Assumption \eqref{eq:ellconv} is satisfied with \eqref{eq:ellmusierv}, 
Theorem \ref{thm:sumBilin} proves ISS of 
\eqref{eq:sirvvaccir} with 
$N(t):=S(t)+I(t)+R(t)$, which 
establishes \eqref{eq:NdissipISSUB} and 
\eqref{eq:NdissipISSUBt} via 
\eqref{eq:NdissipISS}. 
The reminder of the analysis is the same as that of 
the model \eqref{eq:sirvvacci} except that 
the convergence of $I$ to zero does not imply 
$\lim_{t\to\infty}R(t)=0$. Since the scalar $R$-system is ISS, 
it is clear that 
\begin{align}
{B}(1-P)<\mu H
\ \Rightarrow\ \lim_{t\to\infty}x(t)=\left[\frac{B(1-P)}{\mu},0,\frac{BP}{\mu}
\right]^T 
\label{eq:sirv_epidemicequivaccirA}
\end{align}
for any constant $B\in\Rset_+$ and $P\in(0,1)$. 

%%%%%%%%%%

The same modifications to other disease models in the previous sections 
can be possible for modeling the newborn vaccination \cite{DLsencontrdiseas12}. 
Their analysis goes in essentially the same way as the one described above for the SIR model. 

%%%%%%%%%%

If non-newborns/non-immigrants are vaccinated \cite{OGR_IDGoptivacci02,ZAM_IDSIRvacci08}, 
a way to modify the SIR model is 
\begin{subequations}\label{eq:sirvvaccis}
\begin{align}
\dot{S}=&B -\rho S -\mu S -\beta IS
\label{eq:sirvvaccisS}
\\
\dot{I}=&\beta IS -\gamma I-\mu I
\label{eq:sirvvaccisI}
\\
\dot{R}=&\gamma I-\mu R
\label{eq:sirvvaccisR}
\\
\dot{A}=&\rho S-\mu A, 
\label{eq:sirvvaccisA}
\end{align}
\end{subequations}
where the constant $\rho\in\Rset_+$ is the vaccination rate. 
The analysis is the same as that of \eqref{eq:sirvvacci} except that 
ISS of S-system with respect to the input $B$ yields 
property \eqref{eq:propxUISS} with $V_U(x_U)=S$, $x_U=S$ 
for $\alpha_U(s)=(\rho+\mu) s$ and $\sigma_U(s)=s$. Hence, 
\begin{align}
{B}<(\rho+\mu) H
\ \Rightarrow\ \lim_{t\to\infty}x(t)=\left[\frac{B(1\!-\!P)}{\mu},0,0,\frac{BP}{\mu}\right]^T 
\label{eq:sirv_epidemicequivaccisA}
\end{align}
for any constant $B\in\Rset_+$ and $\rho\in\Rset_+$. 
Thus, the disease can be eradicated by a sufficiently high vaccination rate $\rho$. 
Irrespective of ${B}<(\rho+\mu) H$, 
the ultimate bounds \eqref{eq:NdissipISSUB} and 
\eqref{eq:NdissipISSUBt} hold, and  
the model \eqref{eq:sirvvaccis} also has 
the bifurcation point at $S=H$ given in \eqref{eq:defRsirv}. 

\section{Concluding Remarks}\label{sec:conc}

This paper has investigated popular models of infectious diseases 
from the viewpoint of iISS and ISS. 
It has been shown that behavior of all the models can be 
analyzed uniformly in terms of a asymptotic gain property of ISS, 
and a Strongly iISS component which is not globally ISS, but admits 
a zero local ISS-gain function. 
The outbreak is caused by the Strongly iISS component which 
is not ISS. However, the disease is eradicable since 
the Strongly iISS component possesses zero local ISS-gain which 
takes effect if a characteristic value is below a threshold. 
The notions of (i)ISS absorb changes of equilibria, and 
provide a module-based framework. 
The analysis of global properties does 
not require direct and heuristic construction of different Lyapunov 
functions of the entire network depending on equilibria. 
This demonstration is the main contribution of this paper. 
The same procedure and explanation are valid even in the presence of 
an outer-loop caused by short-time immunity. 
The source of the particular iISS component is bilinearity. 
Indeed, scalar linear systems can never exhibit peaks, the outbreak. 
Although the bilinearity is the only nonlinearity in the popular simplest 
models, the theoretical tools presented in this paper accommodate a broad 
class of nonlinearities, such as saturation,
non-monotone nonlinearities 
\cite{CAPSERgennonID78,LIULEVgennonID86,KOROgennonID06,XIARUAgennonID07,EnaNakIDlyapdelay11,EnaNakIDlyap14} and others, 
as long as component models retain appropriate iISS and ISS properties. 
In fact, the arguments in this paper 
rely on neither linearity nor particular nonlinearities. 
Only ISS, iISS and ISS-gain characterizations 
are utilized.

This paper has not reported new epidemiologic discoveries. 
Nevertheless, the system and signal treatment is expected to be superior to 
heuristic approaches in finding control strategies  
for eradicating or containing the spread of diseases. 
The proposed option aims to facilitate the research on control design with global guarantees. 
It is worth noticing that Lyapunov functions used in this paper 
are weighted sum of populations, 
which are simpler than logarithmic functions that have been 
popular in the field of mathematical epidemiology \cite{KOROLyap02}. 
More importantly, (i)ISS Lyapunov functions constructed in this paper 
are different from Lyapunov functions in the conventional concept, and 
the construction of Lyapunov functions does not need 
preprocessing of equilibria. 
The vaccination discussed in this paper is open-loop. 
Interesting future research includes introduction of the (i)ISS framework 
to closed-loop control design 
(see, e.g., \cite{NIEfcimpulseID12,DLsencontrdiseas12,ALQoutlinzdiseas12} and references therein). 

%\begin{acknowledgements}
%\end{acknowledgements}

%
\end{document}